\documentstyle[12pt]{article}
\include{mssymb}
\newtheorem{prop}{Proposition}
\newtheorem{thm}{Theorem}
\newtheorem{cor}{Corollary}
\newtheorem{lemma}{Lemma}

\def\endproof{$\triangle$}
\def\P{{\Bbb P}}

\def\G{{\Bbb G}}
\def\Z{{\Bbb Z}}
\def\A{{\Bbb A}}

\def\Q{{\Bbb Q}}

\begin{document}
\title{Characteristic classes in the Chow ring}
\author{Dan Edidin and William Graham\thanks{The
authors were partially supported by NSF postdoctoral fellowships}
\\ Department of Mathematics\\ University of Chicago\\ Chicago IL
60637}
\date{}
\maketitle
Let $G$ be a reductive algebraic group over an algebraically closed
field $k$.  An algebraic characteristic class of degree $i$ for principal
$G$-bundles on schemes is a function $c$
assigning to each principal $G$-bundle $E \rightarrow X$ an element
$c(E)$ in the Chow group $A^iX$, natural with respect to pullbacks.
These classes are analogous to topological characteristic
classes (which take values in cohomology), and two natural questions
arise.  First, for smooth schemes there is a natural map from the Chow
ring to cohomology, and we can ask if topological characteristic
classes are algebraic.  Second, because the
notion of algebraic principal bundles on schemes is more restrictive
than the notion of topological principal bundles, we can ask if there
are algebraic characteristic classes which do not come from
topological ones.  For example, for rank $n$ vector bundles
(corresponding to $G=GL(n)$), the only topological characteristic
classes are polynomials in the Chern classes, which are represented by
algebraic cycles, but until now it was not known if these were the
only algebraic characteristic classes (\cite[Problem (2.4)]{Vistoli}).
In this paper we describe the ring of algebraic characteristic classes
and answer these questions.

One subtlety which does not occur in topology is that there are
two natural notions of algebraic principal
$G$-bundles on schemes, those which are locally trivial in the \'etale
topology and those which are locally trivial in the Zariski topology.
(Of course for groups which are special in the sense of
\cite{Sem-Chev}, all principal bundles are Zariski locally trivial.
Tori, $GL(n)$, $SL(n)$, and $Sp(2n)$ are all examples of special
groups.)
Let ${\cal C}(G)$ denote the ring of characteristic classes for principal
$G$-bundles locally trivial in the \'etale topology,
and ${\cal C}_{Zar}(G)$ the analogous ring for bundles locally
trivial in the Zariski topology. Since any bundle which is locally
trivial in the Zariski topology is locally trivial in the \'etale
topology, there is a natural homomorphism ${\cal C}(G) \rightarrow
{\cal C}_{Zar}(G)$, but the two rings can differ, as we discuss below.

To state our results, we need some notation.  Let $T$
be a maximal torus of $G$, $B$ a Borel subgroup containing $T$, $U$
the unipotent radical of $B$, and $W$ the
Weyl group.  Let $\hat{T}$ denote the group of
characters of $T$. Let $S(\hat{T})$
be the symmetric algebra over $\Z$ of the free abelian group
$\hat{T}$.
If $\lambda : T \rightarrow \G_m$
is a character of $T$, we extend $\lambda$ to a character of $B$
(still denoted $\lambda$) by making $\lambda |_{U}$ trivial.  Let
$k_{\lambda}$ be the corresponding 1-dimensional representation of $B$
on the ground field $k$.  If $E
\rightarrow X$ is a principal $G$-bundle, we can form a line bundle
$L_{\lambda} = E \times^{B} k_{\lambda} \rightarrow E/B$.   The map
${\lambda} \mapsto c_{1}(L_{\lambda})$ is a group
homomorphism from $\hat{T}$ to $A^{1}(E/B)$.
This extends to a ring homomorphism $ \Phi_E: S(\hat{T}) \rightarrow
A^{*}(E/B)$.

\begin{thm}\label{t.bgzar}
$(a)$ Let $E \rightarrow X$ be a Zariski locally trivial principal
$G$-bundle, and let $f \in S(\hat{T})^W$.  There is a unique element $c_{f}(E)
\in
A^{*}X$  which pulls back to $ \Phi_E(f) \in A^{*}(E/B).$
The
assignment $E \mapsto c_{f}(E)$ is a characteristic class for Zariski
locally trivial principal $G$-bundles.  If $E \rightarrow X$ is not
Zariski locally trivial, the same statements hold after tensoring with $\Q$.

$(b)$ The map $S(\hat{T})^W \rightarrow {\cal C}_{Zar}(G)$,
$f \mapsto c_{f}$, is
an isomorphism.

$(c)$ The map $S(\hat{T})^W \otimes \Q \rightarrow {\cal C}(G) \otimes \Q$,
$f \mapsto c_{f}$, is an isomorphism.

\end{thm}
The following is an immediate consequence of our result for $G=GL(n)$.
\begin{cor} \label{chern.cor}
The only algebraic characteristic classes for vector bundles
are polynomials in the Chern classes.
\end{cor}

In topology the ring of characteristic classes is just $H^*(BG)$
where $BG$ is the classifying space of $G$. Since
$H^*(BG;\Q) = S(\hat{T})^W \otimes \Q$, the theorem says that
${\cal C}(G) \otimes \Q \simeq {\cal C}(G)_{Zar} \otimes \Q \simeq H^*(BG,\Q)$.
Borel showed that with integer coefficients $H^*(BG;\Z)/Tors$
injects into $S(\hat{T})^W$ (cf. \cite{Feshbach}). This implies the
following corollary.
\begin{cor}
All characteristic
classes (mod torsion) from topology are algebraic when applied
to Zariski locally trivial principal $G$-bundles.
\end{cor}

{\bf Discussion.} Vistoli (\cite{Vistoli}, cf.
\cite{Atiyah-Hirzebruch}) constructs an injective map from $S(T)^W
\otimes \Q \rightarrow {\cal C}(G)
\otimes \Q \simeq {\cal C}_{Zar}(G) \otimes \Q$.
What is new in Theorem 1(c) is that these are the only characteristic
classes with rational coefficients.

Given a
representation $V$ of $G$ we can construct characteristic
classes
which assign to a principal $G$-bundle $E \rightarrow X$
a polynomial in the Chern classes of the associated vector bundle $E \times
^{G} V \rightarrow X$. The classes Vistoli constructs are
$\Q$-linear combinations of these. Combining Vistoli's
result with Theorem 1(c) shows that every element of
${\cal C}(G) \otimes \Q$ is of this form.

In general, we do not know if integral polynomials in Chern classes of
associated vector bundles generate ${\cal C}(G)$ (in part because
because we do not have a complete description of ${\cal C}(G)$ for
non-special groups).  However, for the classical groups other than
$SO(2n)$, ${\cal C}_{Zar}(G)$ is generated in this way.  For
$G=SO(2n)$ the only generator of ${\cal C}_{Zar}(G)$ which cannot be
so constructed is the Euler class, which was constructed in \cite{E-G}
(this paper gives an implicit construction of the Euler class as
well).  For the exceptional groups, we do not know which
characteristic classes come from representations.

In topology, the ring of characteristic classes with arbitrary
coefficients can be computed from the ring of characteristic classes
with integer coefficients.  By contrast, we do not know how to use
Theorem 1 to compute the ring ${\cal C}_{Zar}(G;R)$ of characteristic
classes with values in $A^*(\mbox{ };R)$, for an arbitrary ring $R$.
Note that ${\cal C}_{Zar}(G;R)$ need not be isomorphic to ${\cal
C}_{Zar}(G) \otimes R$.  $G=SO(n)$ is an example, since ${\cal
C}_{Zar}(G;\Z /2 \Z )$ contains Stiefel-Whitney classes (\cite{E-G}),
which are not in ${\cal C}_{Zar}(G) \otimes \Z /2 \Z $.

{\bf Examples.} For $G$ equal to $GL(n)$, $Sp(2n)$, or $SO(2n+1)$,
${\cal C}_{Zar}(G)$ is isomorphic to the polynomial ring
$\Z[t_1,\ldots,t_n]$, where the $t_i$ are constructed from the
defining representation of $G$, as follows.  For $GL(n)$, we can take
$t_i$ to be the $i$-th Chern class of the associated vector bundle.
For $Sp(2n)$ and $SO(2n+1)$, $t_i$ is the $2i$-th Chern class.

For $G = SO(2n)$, ${\cal C}_{Zar}(G)$ is also isomorphic to the
polynomial ring $\Z[t_1,\ldots,t_n]$.  For $i < n$, $t_i$ is the
$2i$-th Chern class of the associated rank $2n$ vector bundle, but
$t_n$ is the Euler class; this class
satisfies the relation $t_{n}^{2} = (-1)^n c_{2n}$.  The Euler class
cannot be expressed as a polynomial in the Chern classes of an
associated vector bundle.  If it could, then it would be defined in
${\cal C}(G)$.  However, Totaro \cite{T2} has shown that in general
only some multiple of the Euler class is defined in ${\cal C}(G)$.
This example shows that the rings ${\cal C}(G)$ and ${\cal
C}_{Zar}(G)$ are not isomorphic in general.

{}From the examples above, one might suspect that ${\cal C}_{Zar}(G)$ is
always a polynomial ring. However, Feshbach \cite{Feshbach} shows that
for $G=$ Spin(10), $H^*(BG)/Tors$ is not a polynomial ring, but is
isomorphic to $S(\hat{T})^W$. Since we prove ${\cal C}_{Zar}(G) =
S(\hat{T})^W$, it follows that ${\cal C}_{Zar}(G)$ is not always
a polynomial ring.  Feshbach also shows that the map $H^*(BG)/Tors
\rightarrow S(\hat{T})^W$ is not surjective for $G$ equal to Spin(11)
or Spin(12).  This shows that there are characteristic classes for
Zariski locally trivial $G$-bundles which cannot be defined in topology.

\medskip

{\bf Remark on the proof.} We only present the proof of Theorem 1
in the Zariski locally trivial case. The same proof works
for general principal bundles after tensoring with $\Q$. The
main point is that if $f: Y \rightarrow X$ is a Zariski
locally trivial flag bundle then $A^*Y$ is free over
$A^*X$. However, for general flag bundles this is not true.

\medskip

{\bf Notation:} A ``scheme" will mean an algebraic scheme over a fixed
algebraically closed field $k$. $A^*X$ refers to the operational Chow
ring of \cite[Chapter 17]{Fulton}. For smooth schemes, $A^*X$
coincides with the usual Chow ring.  $A_*X$ refers to the Chow groups
of \cite[Chapter 2]{Fulton}.  The evaluation of a characteristic class
$c$ on a principal bundle $E \rightarrow X$ will be denoted $c(E
\rightarrow X)$, or simply $c(E)$ when there is no confusion about the
base; $c(E \rightarrow X)$ is always an element of $A^*X$.

\medskip

{\bf Acknowledgements:} We would like to thank Burt Totaro for
generously sharing his ideas with us, without which this paper would
have been impossible.  The authors also benefitted from conversations
with Rahul Pandharipande.

\paragraph{Reduction to the quasi-projective case} We can define the
ring of characteristic classes for principal $G$-bundles over
quasi-projective schemes as above.  Note that this ring might not be
the same as the ring of characteristic classes of bundles over all
schemes, since there could be characteristic classes over
quasi-projective schemes that cannot be defined for all schemes.
Conversely, there might be nonzero characteristic classes which are
zero when applied to any bundle over a quasiprojective scheme.  The
following lemma says that neither of these possibilities occurs.

\begin{lemma}
Let $G$ be an algebraic group.
The ring of characteristic classes of principal $G$-bundles over
quasi-projective schemes is the same as the ring of characteristic
classes of principal $G$-bundles over arbitrary schemes.  The same
statement is true for characteristic classes of Zariski locally
trivial principal $G$-bundles.
\end{lemma}

Recall (\cite[Definition 18.3 and
Lemma 18.3]{Fulton}) that every scheme $X$ has a {\it Chow envelope}
$p:X^\prime \rightarrow X$, which is a quasi-projective scheme such
that for every subvariety $V \subset X$ there is a subvariety
$V^\prime \subset X^\prime$ mapping birationally to $V$.  The most
important
property of Chow envelopes is that the pullback $p^*:A^*X \rightarrow
A^*X^\prime$ is injective.
Moreover,
$X^\prime$ may be chosen so that there exists an open set $U \subset
X$ such that $p$ is an isomorphism over $U$ and the complement $Z = X
- U$ has smaller dimension than $X$.  For the sake of brevity, we will
call such an $X'$ a {\it birational} Chow envelope of $X$. In this
case, the image of $p^*$ in $A^*X'$ has been described by
Kimura (\cite{Kimura}).

Proof: The proof of the lemma is in two steps.  We first show that if a
characteristic class $c$ for bundles over arbitrary schemes is zero
for all bundles over quasi-projective schemes, then it is zero.  We
then show that a characteristic class $c$ for bundles over
quasi-projective schemes is the restriction to such bundles of a
characteristic class for bundles over arbitrary schemes.

Step 1: Suppose $E \rightarrow X$ is a principal bundle over an
arbitrary scheme.  Let $p:X^\prime \rightarrow X$ be a Chow envelope,
and let $E^\prime \rightarrow X^\prime$ be the pullback bundle.  By
hypothesis, $0 = c(E' \rightarrow X') = p^*c(E \rightarrow X) $.  Since
the pullback $p^*:A^*X \rightarrow A^*X^\prime$ is injective, $c(E
\rightarrow X) = 0$ as desired.

Step 2: It suffices to show the following.  Suppose $E \rightarrow X$
is a $G$-bundle and $p:X' \rightarrow X$ a map with $X'$ quasiprojective.
Let $c$ be a characteristic
class for bundles over quasi-projective schemes.  Then $c(p^*E
\rightarrow X') \in A^*X'$ is the pullback of a class in $A^*X$.

Suppose the statement is true for all schemes of dimension smaller
than $X$.  We will first prove the statement assuming that $X'
\rightarrow X$ is a birational Chow envelope.  Let $Z' =
p^{-1}Z$, and consider the following diagram:
$$\begin{array}{ccc}
Z'& \stackrel{i'}\rightarrow & X'\\
\scriptsize{q} \downarrow & & \scriptsize{p} \downarrow\\
Z & \stackrel{i} \rightarrow X
\end{array}$$

By \cite[Theorem 3.1]{Kimura}, to show that $c(p^*E \rightarrow X')
\in A^*X'$ is the pullback of a class from $A^*X$ is equivalent to
showing that $i'^*c(p^*E \rightarrow X')$ is the pullback of a class
from $A^*Z$.  But $i'^*c(p^*E \rightarrow X') = c(q^*i^*E \rightarrow
Z')$, which is a pullback from $A^*Z$ by the inductive hypothesis.
This proves the statement if $X' \rightarrow X$ is a birational Chow
envelope. To deduce the statement for arbitrary $X'$, let $\tilde{X}
\rightarrow X$ be a birational Chow envelope, and consider the Cartesian
diagram: $$
\begin{array}{ccc}
\tilde{X'}& \stackrel{q'} \rightarrow & X'\\
\scriptsize{\tilde{p}} \downarrow & & \scriptsize{p} \downarrow\\
\tilde{X}&  \stackrel{q} \rightarrow & X
\end{array}$$
Let $E'$, $\tilde{E}$, and $\tilde{E}'$ be the pullback bundles.  Then
there exists $\alpha \in A^*X$ such that $c(\tilde{E}') = \tilde{p}^*q^*
\alpha = q'^*p^* \alpha$.  But $c(\tilde{E}') = q'^*c(E')$.  Since $\tilde{X}'
\rightarrow X'$ is an envelope, $q'^*$ is injective, so $c(E') = p^*
\alpha$, as desired. \endproof

By the above lemma, to calculate rings of characteristic classes it
suffices to do so for bundles over quasi-projective schemes.  For the
remainder of this paper we will therefore assume that all schemes are
quasi-projective.

\paragraph{Classifying Spaces in Algebraic Geometry} The purpose of
this section is to describe a construction, due to Totaro, of an
algebro-geometric substitute for the classifying spaces of topology.
In order to say what we mean by a substitute, it is helpful to recall
what classifying spaces are.  In topology, given a contractible space
$EG$ on which $G$ acts freely, the quotient space $BG$ is called the
classifying space for $G$, and the bundle $EG \rightarrow BG$ is
called the universal bundle.  A principal $G$-bundle over any space
$X$ is pulled back from the universal bundle by a classifying map $X
\rightarrow BG$, and homotopic maps induce the same pullback bundle.
Hence there is an isomorphism of $H^*BG$ with (cohomology)
characteristic classes of principal $G$-bundles.

Because the topological spaces $EG$ and $BG$ are infinite dimensional,
it seems unreasonable to expect to find a bundle of schemes which is
universal in the category of schemes.  The next best thing would be to
find a directed system of bundles of schemes $E_{n} \rightarrow
B_{n}$, such that any principal $G$-bundle of schemes pulls back from
some $E_{n} \rightarrow B_{n}$.  However, we do not know how to
produce this.\footnote{There is a classifying space $BG_{\cdot}$ which
is a simplicial scheme whose Chow groups can be defined by the work of
Bloch \cite{Bloch}. However, computing these Chow groups seems
difficult.} What Totaro produces instead is a directed system of
bundles as above, with the following property: for any principal
$G$-bundle $E \rightarrow X$, there is a map $X' \rightarrow X$ with
fibers isomorphic to $\A ^{m}$ such that the pullback bundle $E'
\rightarrow X'$ is pulled back from one of the bundles in the directed
system by a map $X' \rightarrow B_{n}$.  Moreover, in this situation,
for fixed $i$, the Chow groups $A^i (B_{n})$ are all isomorphic for
$n$ sufficiently large.  By abuse of notation, we call this group
$A^i(BG)$.

The directed system above is Totaro's algebro-geometric substitute for
$BG$.  The system possesses the following two properties, analogous to
ones familiar from topology.  First, there is a map $A^*(BG)
\rightarrow {\cal C}(G)$, $\alpha \mapsto c_{\alpha}$, defined as
follows.  Because the bundle $X' \rightarrow X$ has fibers isomorphic
to $\A ^{m}$, by \cite[Theorem 8.3]{Gillet}, the Chow groups $A^*X$
and $A^*X'$ are isomorphic.  Thus, given $\alpha \in A^i(BG)$,
pullback gives a class, which we denote $c_{\alpha}(E)$ in $A^{i}X'
\cong A^{i}X$.  The class $c_{\alpha}(E)$ in $A^{i}X$ depends only on
the class in $A^i(BG)$ and on the bundle $E \rightarrow X$, not on the
choice of $X' \rightarrow X$ or the classifying map $X' \rightarrow
B_{n}$.  Moreover, the assignment $E \mapsto c_{\alpha}(E)$ is natural
with respect to pullbacks.  Hence $c_{\alpha}$ is a characteristic
class.  Now arguing as in topology shows the following.
\begin{thm}(Totaro)
The map $A^*(BG)
\rightarrow {\cal C}(G)$ is an isomorphism.
\end{thm}
Because of this theorem, we will use the notations $A^*(BG)$ and
${\cal C}(G)$ interchangeably.

The other nice property that this directed system possesses is simply
that by replacing $X$ by $X'$ one can assume that the bundle $E
\rightarrow X$ is pulled back from some bundle in the directed system.
This enables one to reduce some computations to computations in the
directed system.

By abuse of notation we will write $EG \rightarrow BG$ as schemes,
where we really mean we are working in some bundle in the directed
system.  This should lead to no confusion, since we will never use the
topological notion of classifying space.

We give some details of Totaro's construction, since we will need it
to compute the groups $A^*(BG)$.  If $V$ is a representation of $G$
let $V^{s}$ be the set of points $v \in V$ such that $G \cdot v$ is
closed in $V$.  Let $V^{sf}$ be the set of $v \in V^{s}$ with trivial
stabilizer.  $V^{sf}$ is an open (possibly empty) subset of $V$.
Since the orbits in $V^{sf}$ are closed in $V$ and have the same
dimension as $G$, there is a geometric quotient $V^{sf} \rightarrow
V^{sf}/G$ (cf. \cite[Chapter 1, Appendix B]{GIT}).  Furthermore, since
the action of $G$ on $V^{sf}$ is free, $V^{sf} \rightarrow V^{sf}/G$
is actually a principal $G$-bundle.  These bundles form the directed
system; the morphisms are just inclusions.  Given a
representation $V$ such that $V^{sf}$ is non-empty (we will see below
that such representations exist), let $W = V \oplus V$.  Then
$(V^{sf} \oplus V) \cup (V \oplus V^{sf}) \subset W^{sf}$.  Hence by
replacing $V$ by a direct sum of copies of $V$ we may assume that $V -
V^{sf}$ has arbitrarily high codimension.  Totaro shows that
$A^{i}(V^{sf}/G)$ does not depend on the choice of $V$ as long as the
codimension of $V - V^{sf}$ is greater than $i$.  His argument
uses Bogomolov's double fibration (cf.
the proof of Lemma \ref{char.tor}).

The classifying map is constructed as follows.  Suppose $E \rightarrow
X$ is a principal $G$-bundle.  Then there exists a bundle $X'
\rightarrow X$, with fibers $\A ^{m}$, such that the total space $ E'
= E \times_X X'$ is affine.  Because $E'$ is affine, there is a
representation $V$ of $G$ and a $G$-equivariant embedding $E'
\rightarrow V$.  ($V$ may be taken to be a subspace of the space of
regular functions on $E'$.)  Since the action of $G$ on $E'$ is closed
and free, the map $E' \rightarrow V$ has image in $V^{sf}$ (in
particular showing there exist $V$ with $V^{sf}$ non-empty).  Thus $E'
\rightarrow X'$ is the pullback of the principal $G$ bundle $V^{sf}
\rightarrow V^{sf}/G$.

Of particular importance is the case when $G=T$ is a torus.  Then
$A^*BT = {\cal C}_{Zar}(T)$, since all principal $T$-bundles in
algebraic geometry are locally trivial in the Zariski topology
\cite{Sem-Chev}.  There is a group homomorphism $\hat{T} \rightarrow
A^1BT$ taking a character of $T$ to the first Chern class of the
associated line bundle on $BT$.  This extends to a ring homomorphism
$S(\hat{T}) \rightarrow A^*BT$.
\begin{lemma}  \label{char.tor}
The map $S(\hat{T}) \rightarrow A^*BT$ is an isomorphism.
\end{lemma}
Proof: Assume first that $T=\G_m$. In this case, we can take $V$ to be
a product of $N$ copies of the 2 dimensional representation acted on
by the torus with weights 1 and $-1$. The complement of $(V^{sf})$ has
codimension $2N-1$.

Claim: For $i < 2N-1$, $A^i(V^{sf}/T) = A^i(\P^{2N-1})$.

Proof of claim: We use Bogomolov's double
fibration argument. Let $W$ be the $2N$-dimensional representation
of $T$ with all weights equal to 1. Then $(W-\{0\})/T=\P^{2N-1}$, and
$(V+(W-\{0\}))/T$ is an affine bundle over $\P^{2N-1}$ so
$A^i((V+(W-\{0\}))/T)=A^i(\P^{2N-1})$. On the other hand
$A^i((V^{sf} +W)/T) =A^i(V^{sf}/T)$. Thus for $i < 2N-1$,
$A^i(V^{sf}/T) = A^i(\P^{2N-1})$.

Hence
$A^*BT = \Z[t]$ where $t$ is a generator of $A^1(\P^n)$.
$S(\hat{T})$ is also a polynomial ring in one variable and the above map
takes a generator to a generator.  This proves the lemma for $T=\G_m$.
In the general case, $T$ will be a product of 1-dimensional tori
$T_1,\ldots,T_n$, and $S(\hat{T})$ decomposes as the tensor product of the
$S(\hat{T}_i)$.  On the other hand, $BT$ is the direct product of the
$BT_i$
and $A^*BT$ decomposes as the tensor product of the
$A^*(BT_i)$.  The map $S(\hat{T}) \rightarrow A^*BT$ is compatible
with these tensor product decompositions and hence is an isomorphism.
\endproof

We will also be interested in the case of $B$-bundles,
where $B$ is a Borel subgroup of $G$. As is the case
with tori, principal $B$-bundles are locally trivial
in the Zariski topology, so we do not have to distinguish between
$A^*BB$ and ${\cal C}_{Zar}(B)$.  If $E \rightarrow X$ is
a $B$-bundle, then $E/U \rightarrow X$ is a $T$-bundle
and $E \rightarrow E/U$ is an affine bundle. In particular,
characteristic classes of $B$-bundles are uniquely determined
by a characteristic class for $T$-bundles. This fact can
be summarized in the following useful lemma.
\begin{lemma}
$A^*BB = A^*BT = S(\hat{T})$. \endproof
\end{lemma}

\paragraph{A bound on characteristic classes of arbitrary
bundles}

Define a map from principal $T$-bundles to Zariski locally
trivial principal $G$-bundles by associating to a $T$-bundle
$E \rightarrow X$ the $G$-bundle $E \times^T G \rightarrow X$.
The associated bundle is Zariski locally trivial, since the original
$T$-bundle is Zariski locally trivial by \cite{Sem-Chev}. Next
define a map ${\cal C}_{Zar}(G) \rightarrow A^*BT$, $c \mapsto c_T$
by the formula
$$c_T(E \rightarrow X)  = c(E \times^T G \rightarrow X).$$

The Weyl group $W$ acts on the set of principal $T$-bundles as
follows: If $Y \rightarrow X$ is a $T$-bundle define $wY \rightarrow
X$ as having total space $Y$ but with twisted $T$ action given by $$y
\cdot_w t = y(wtw^{-1})$$ (where the action on the right is the
original action).  The formula $wc(Y \rightarrow X) = c(w^{-1}Y
\rightarrow X)$ defines a $W$ action on $A^*BT$ which is compatible
with the $W$ action on the representation ring $S(\hat{T})$.
\begin{lemma}\label{l.inj}
The map ${\cal C}_{Zar}(G) \rightarrow A^*BT$, $c \mapsto c_T$ is
injective, with image in $A^*BT^W$.
\end{lemma}
Proof: Suppose $c_T=0$. Let
$p:E \rightarrow X$ by any Zariski locally trivial $G$-bundle.  Let
$q:E/T \rightarrow X$ be the associated $G/T$-bundle.  Now, $E
\rightarrow E/T$ is a $T$-bundle.  By assumption, $c_T(E \rightarrow E/T) = 0$.
Thus, by definition, $c(E \times^T G \rightarrow E/T) = 0$. But $E
\times^T G \simeq q^*E$ as $G$-bundles over $E/T$. Hence $c(q^*E
\rightarrow E/T)= q^*(c(E \rightarrow T)) = 0$.  But $q^*$ is
injective because $q$ factors into a composition $E/T \rightarrow E/B
\rightarrow X$; the first map is an affine bundle and the second map
is a proper locally trivial bundle.  Thus $c(E \rightarrow X) = 0$.
Since $E \rightarrow X$ was arbitrary, $c=0$, so the map is injective.
To see that $c_T$ is $W$-invariant notice that the map $$(w^{-1}Y)
\times^T G \rightarrow Y \times^T G$$ sending $(y,g) \mapsto
(y,w^{-1}g)$ is an isomorphism of $G$-bundles.  Hence $wc_T = c_T$
\endproof

\paragraph{Proof of the theorem}
By Lemma \ref{l.inj} we know that
${\cal C}_{Zar}(G) \subset S(\hat{T})^W$. It remains
to show that any $W$-invariant polynomial $f \in S(\hat{T})$
defines a characteristic class. In particular
we must show that if $E \rightarrow X$ is a locally trivial
$G$-bundle then the class $p=\Phi_E(f)$ in $A^*(E/B)$ is
a pullback from $X$.

By \cite{Kimura}, because the $G/B$-bundle $E/B \rightarrow X$ is
locally trivial, showing
$\Phi_E(f)$ is a pullback is equivalent
to showing that $\pi_1^*p=\pi_2^*p$ in the following diagram.
$$\begin{array}{ccc}
E/B \times_X E/B & \stackrel{\pi_1} \rightarrow & E/B\\
\downarrow \scriptsize{\pi_2} & & \downarrow\\
E/B & \rightarrow & X
\end{array}$$
The theorem follows from the following lemma.
\begin{lemma} \label{l:trick}
$\pi_1^* p = \pi_{2}^* p.$
\end{lemma}
Proof: First note that by \cite{Vistoli} (cf. \cite{Atiyah-Hirzebruch}),
$S(\hat{T})^W \otimes \Q \subset A^*BG \otimes \Q \subset
{\cal C}_{Zar}(G) \otimes \Q$. Thus, there is some
constant $n \in \Z$ such that $nf$ is a characteristic
class, so $n(\pi_1^*p -\pi_2^*p)=0$ for any principal bundle
$E \rightarrow X$ (this bundle need not be locally trivial).

Since the bundle $E \rightarrow X$ is pulled back from the
universal bundle over $BG$, it suffices to prove the lemma for
$X=BG$. In this case, $E/B = BB$ and our fiber square becomes
$$\begin{array}{ccc}
BB \times_{BG} BB & \stackrel{\pi_1} \rightarrow & BB\\
\downarrow \scriptsize{\pi_2} & & \downarrow\\
BB & \rightarrow & BG
\end{array}$$
Now, $BB \times_{BG} BB
\rightarrow BB$ is a fiber bundle with fiber $G/B$.  The structure
group of this bundle reduces to $B$, so this bundle is Zariski locally
trivial.  Since $G/B$ has a decomposition into affine cells, Proposition
\ref{l.loctriv} below implies that
$A^{*}(BB \times_{BG} BB) \cong
A^{*}(G/T) \otimes A^{*}BB$.  Since $A^{*}(G/B)$ and $A^{*}BB$
are both torsion-free, so is $A^{*}(BB \times_{BG} BB)$.  Thus,
$\pi_{1}^* p - \pi_{2}^* p=0$, as desired. \endproof

\begin{prop} \label{l.loctriv}
Let $f: Y \rightarrow X$ be a smooth proper locally trivial fibration with $X$
smooth,
whose fiber $F$ has a decomposition into affine cells. Then $A^*Y$ is
(non-canonically) isomorphic to $A^*X \otimes A^*F$ as an abelian group.
\end{prop}
Proof:
Since everything is
smooth we may identify the Chow cohomology $A^*Y$ with the Chow
homology $A_*Y$.
Let $U \subset X$ be an open set over which $f$ is trivial. Since
$F$ has a cellular decomposition, $A^*(f^{-1}(U)) \simeq A^*U \otimes A^*F$
(cf. \cite[Example 1.10.2]{Fulton}).
Let
$b_1 ,  \ldots , b_n$ be a basis for
 $A_*F_{x_0}$ where $F_{x_0}$ is the fiber over a point
$x_0 \in U$. Choose representatives
$b_i = \sum a_{ij} [V_{ij}]$, and set $B_i = \sum a_{ij} [\overline{U \times
V_{ij}}] \in A_*(Y)$. We claim the $B_1, \ldots , B_n$ form a basis
for $A_*Y$ over $A_*X$.

To prove the claim argue as follows. Let $F_x$ be any fiber. Since
$F_x$ is algebraically equivalent to any other fiber, the restrictions of
the $B_i$ to $F_x$ are algebraically equivalent to their restrictions
to $F_{x_{0}}$. However, because $F$ has an affine cellular decomposition, the
group of cycles modulo algebraic equivalence is equal to the Chow group.
Thus, the classes $B_i$ restrict to a basis for the Chow group of every fiber.

The proposition now follows from the following algebraic version of
the Leray-Hirsch theorem.
\begin{lemma} \label{l.lh}
Let $f: Y \rightarrow X$ be a smooth proper locally trivial fibration
whose fiber $F$ has a decomposition into affine cells.
Let $\{B_i\} \in A^*Y$ be a collection of classes that restrict to a basis
of the Chow groups of the fibers. Then $\{B_i\}$ form a basis for
$A_*Y$ over $A_*X$, i.e. every $y \in A_*Y$ has a unique expression of
the form $y= \sum B_i \cap f^*x_i$ with $x_i \in A_*X$ (This makes sense
even when $X$ is singular, so that $A_*X$ is not a ring.)
\end{lemma}

Proof: This is an easy generalization of \cite[Lemma 2.8]{E-S}.
\endproof

\end{document}